\documentclass{article} 

\usepackage{amsmath, amssymb}

\newcommand{\RM}{\mathbb{R}}

\title{METZLER/ZETA CORRESPONDENCE}

\author{Yusuke IDE \\ 
Department of Mathematics, College of Humanities and Sciences, \\
Nihon University \\ 
Sakura-josui, Setagaya-ku, Tokyo 156-8550, JAPAN \\ 
\\
Takashi KOMATSU \\
Math. Research Institute Calc for Industry, \\
Minami, Hiroshima, 732-0816, JAPAN \\ 
e-mail: ta.komatsu@sunmath-calc.co.jp \\ 
\\
Norio KONNO \\
Department of Applied Mathematics, Faculty of Engineering, \\ 
Yokohama National University \\
Hodogaya, Yokohama 240-8501, JAPAN \\
e-mail: konno-norio-bt@ynu.ac.jp \\ 
\\
Iwao SATO \\ 
Oyama National College of Technology \\
Oyama, Tochigi 323-0806, JAPAN \\ 
e-mail: isato@oyama-ct.ac.jp }
 
\begin{document}
 \maketitle

\clearpage

\begin{abstract} 
We present an explicit formula for the determinant on the Metzler matrix of a digraph $D$. 
Furthermore, we introduce a walk-type zeta function with respect to this Metzler matrix of the symmetric digraph 
of a finite torus, and express its limit formula by using the integral expression. 
\end{abstract}

\vspace{5mm}


\vspace{5mm}

The contact author for correspondence:

Iwao Sato 

Oyama National College of Technology, 
Oyama, Tochigi 323-0806, JAPAN 

E-mail: isato@oyama-ct.ac.jp

\clearpage

\section{Introduction}

Ihara \cite{I} defined the Ihara zeta functions of graphs, and showed that the 
reciprocals of the Ihara zeta functions of regular graphs are explicit polynomials. 
The Ihara zeta function of a regular graph $G$ associated with a unitary 
representation of the fundamental group of $G$ was developed 
by Sunada \cite{Sun1, Sun2}. 
Hashimoto \cite{H} treated multivariable zeta functions of bipartite graphs. 
Furthermore, Hashimoto \cite{H} gave a determinant expression for the Ihara zeta function 
of a general graph by using the edge matrix. 
Bass \cite{Bass} generalized Ihara's result on the Ihara zeta function of 
a regular graph to an irregular graph $G$. 
Stark and Terras \cite{ST} gave an elementary proof of Bass' theorem, and 
discussed three different zeta functions of any graph. 
Various proofs of Bass' theorem were given by Kotani and Sunada \cite{KS},   
and Foata and Zeilberger \cite{FZ}.

Recently, there were exciting developments between quantum walk \cite{Ambainis2003, Kempe2003, 
Kendon2007, Konno2008b, VA} 
on a graph and the Ihara zeta function of a graph. 
We investigated a new class of zeta functions for many kinds of walks including the quantum walk (QW) and 
the random walk (RW) on a graph by a series of ''Zeta Correspondence" of our previous work \cite{K1, K2, K3, K4, K5, K6, K7}. 
In Walk/Zeta Correspondence \cite{K2}, 
a walk-type zeta function was defined without use of the determinant expressions of zeta function 
of a graph $G$, and various properties of walk-type zeta functions of RW, correlated random walk (CRW) 
and QW on $G$ were studied. 
Also, their limit formulas by using integral expressions were presented.

The Ihara zeta function has many applications in both pure and applied mathematics, 
including for instance dynamical systems, spectral graph theory and complex network analysis. 
In complex network analysis, Masuda, Preciado and Ogura \cite{MPO} studied the stochastic 
susceptible-infected-susceptible model (SIS model) of epidemic processes on finite directed and weighted networks, 
and presented a new lower bound on the exponential rate by using the leading eigenvalue of some matrix 
composed from the digraph version of the edge matrix. 
The edge matrix appeared in the determinant expressions for the Ihara zeta function. 
Furthermore, the above matrix is a kind of the Metzler matrix. 

In this paper, we present an explicit formula for a determinant on the Metzler matrix of a digraph $D$. 
Moreover, we introduce a walk-type zeta function with respect to this Metzler matrix of the symmetric digraph 
of a finite torus, and express its limit formula by using integral expression. 

This paper is organized as follows:
In Section 2, we give a short review about the Ihara zeta function of a graph. 
In Section 3, we review "Walk/Zeta Correspondence" studied in \cite{K2}. 
In Section 4, we state a new lower bound on the exponential rate by Masuda, Preciado and Ogura \cite{MPO}. 
In Section 5, we present an explicit formula for the determinant on the above Metzler matrix of a digraph $D$. 
In Section 6, we give an explicit formula for the characteristic polynomial of the above Metzler matrix 
of the symmetric digraph of a graph, and obtain its spectrum. 
In Section 7, we introduce a walk-type zeta function with respect to this Metzler matrix of the symmetric digraph 
of a finite torus, and express its limit formula by using integral expression.

\section{The Ihara zeta function of a graph} 

Graphs and digraphs treated here are finite.
Let $G$ be a connected graph and $D(G)= \{ (u,v),(v,u) \mid uv \in E(G) \} $ 
the arc set of the symmetric digraph 
corresponding to $G$. 
For $e=(u,v) \in D(G)$, set $u=o(e)$ and $v=t(e)$. 
Furthermore, let $e^{-1}=(v,u)$ be the {\em inverse} of $e=(u,v)$. 

A {\em path $P$ of length $n$} in $G$ is a sequence 
$P=(e_1, \cdots ,e_n )$ of $n$ arcs such that $e_i \in D(G)$,
$t( e_i )=o( e_{i+1} )(1 \leq i \leq n-1)$. 
If $e_i =( v_{i-1} , v_i )$ for $i=1, \cdots , n$, then we write 
$P=(v_0, v_1, \cdots ,v_{n-1}, v_n )$. 
Set $ \mid P \mid =n$, $o(P)=o( e_1 )$ and $t(P)=t( e_n )$. 
Also, $P$ is called an {\em $(o(P),t(P))$-path}. 
We say that a path $P=( e_1 , \cdots , e_n )$ has a {\em backtracking} 
if $ e^{-1}_{i+1} =e_i $ for some $i \ (1 \leq i \leq n-1)$. 
A $(v, w)$-path is called a {\em $v$-cycle} 
(or {\em $v$-closed path}) if $v=w$. 
The {\em inverse cycle} of a cycle 
$C=( e_1, \cdots ,e_n )$ is the cycle 
$C^{-1} =( e^{-1}_n, \cdots ,e^{-1}_1 )$.

We introduce an equivalence relation between cycles. 
Two cycles $C_1 =(e_1, \cdots ,e_m )$ and 
$C_2 =(f_1, \cdots ,f_m )$ are called {\em equivalent} if 
$f_j =e_{j+k} $ for all $j$. 
The inverse cycle of $C$ is not equivalent to $C$ if $\mid C \mid \geq 3$. 
Let $[C]$ be the equivalence class which contains a cycle $C$. 
Let $B^r$ be the cycle obtained by going $r$ times around a cycle $B$. 
Such a cycle is called a {\em multiple} of $B$. 
A cycle $C$ is {\em reduced} if 
both $C$ and $C^2 $ have no backtracking. 
Furthermore, a cycle $C$ is {\em prime} if it is not a multiple of 
a strictly smaller cycle. 
Note that each equivalence class of prime, reduced cycles of a graph $G$ 
corresponds to a unique conjugacy class of 
the fundamental group $ \pi {}_1 (G,v)$ of $G$ at a vertex $v$ of $G$. 

The {\em Ihara(-Selberg) zeta function} of $G$ is defined by  
\[
{\bf Z} (G,u)= \prod_{[C]} (1- u^{ \mid C \mid } )^{-1} , 
\]
where $[C]$ runs over all equivalence classes of prime, reduced cycles 
of $G$. 

Let $G$ be a connected graph with $N$ vertices and $M$ edges. 
Then two $2M \times 2M$ matrices 
${\bf B} = {\bf B} (G)=( {\bf B}_{ef} )_{e,f \in D(G)} $ and 
${\bf J}_0 ={\bf J}_0 (G) =( {\bf J}_{ef} )_{e,f \in D(G)} $ 
are defined as follows: 
\[
{\bf B}_{ef} =\left\{
\begin{array}{ll}
1 & \mbox{if $t(e)=o(f)$, } \\
0 & \mbox{otherwise, }
\end{array}
\right.
\   
{\bf J}_{ef} =\left\{
\begin{array}{ll}
1 & \mbox{if $f= e^{-1} $, } \\
0 & \mbox{otherwise.}
\end{array}
\right.
\]
Note that 
\[
{\bf J}_0 = {\bf B} \circ {\bf B}^T , 
\]
where ${\bf B}^T $ is the transpose of ${\bf B}$, and the {\em Schur/Hadamard product} ${\bf A} \circ {\bf B} $ 
of two matrices ${\bf A} $ and ${\bf B} $ is defined by 
\[
( {\bf A} \circ {\bf B} )_{ij} = {\bf A}_{ij} \cdot {\bf B}_{ij} . 
\]
The matrix ${\bf B} - {\bf J}_0 $ is called the {\em edge matrix} of $G$.

\newtheorem{theorem}{Theorem}
\begin{theorem}[Ihara; Hashimoto; Bass] 
Let $G$ be a connected graph with $N$ vertices and $M$ edges. 
Then the reciprocal of the Ihara zeta function of $G$ is given by 
\[
{\bf Z} (G,u) {}^{-1} = \det ( {\bf I}_{2M} - u( {\bf B} - {\bf J}_0 )) 
=(1- u^2 )^{M-N} \det( {\bf I}_N  -u {\bf A} (G)+ u^2 ( {\bf D}_G - {\bf I}_N  )) , 
\]
where ${\bf D}_G =( d_{ij} )$ is the diagonal matrix with 
$d_{ii} = \deg {}_G \  v_i \  (V(G)= \{ v_1 , \cdots , v_N \} )$. 
\end{theorem}

\section{Walk/Zeta Correspondence on torus} 

We state the Walk/Zeta Correspondence on a finite torus. 
At first, we give the definition of the $2d$-state discrete-time walk on $T^d_N$. 
The discrete-time walk is defined by using a {\em shift operator} and a {\em coin matrix} which will be mentioned below.

Let $f : V( T^d_N ) \longrightarrow \mathbb{C}^{2d}$.
For $j = 1,2,\ldots,d$ and ${\bf x} =( x_1 , \ldots , x_d ) \in V( T^d_N )$, the shift operator $\tau_j$ is defined by 
\begin{align*}
(\tau_j f)({\bf x}) = f( {\bf x} - {\bf e}_{j}),
\end{align*} 
where $\{ {\bf e}_1,{\bf e}_2,\ldots,{\bf e}_d \}$ denotes the standard basis of $\mathbb{R}^d$.

Let $A=(a_{ij})_{i,j=1,2,\ldots,2d}$ be a $2d \times 2d$ matrix with $a_{ij} \in \mathbb{C}$ for $i,j =1,2,\ldots,2d$. 
We call $A$ the {\em coin matrix}. If $a_{ij} \in [0,1]$ and $\sum_{i=1}^{2d} a_{ij} = 1$ for any $j=1,2, \ldots, 2d$, 
then the walk is a CRW. 
In particular, when $a_{i1} = a_{i2} = \cdots = a_{i 2d}$ for any $i=1,2, \ldots, 2d$, this CRW becomes a RW. If $A$ is unitary, 
then the walk is a QW. So our class of walks contains RWs, CRWs, and QWs as special models.

To describe the evolution of the walk, we decompose the $2d \times 2d$ coin matrix $A$ as
\begin{align*}
A=\sum_{j=1}^{2d} P_{j} A,
\end{align*}
where $P_j$ denotes the orthogonal projection onto the one-dimensional subspace $\mathbb{C}\eta_j$ in $\mathbb{C}^{2d}$. 
Here $\{\eta_1,\eta_2, \ldots, \eta_{2d}\}$ denotes a standard basis on $\mathbb{C}^{2d}$.
 
The discrete-time walk associated with the coin matrix $A$ on $T^d_N$ is determined by the $2d N^d \times 2d N^d$ matrix
\begin{align}
M_A=\sum_{j=1}^d \Big( P_{2j-1} A \tau_{j}^{-1} + P_{2j} A \tau_{j} \Big).
\label{unitaryop1}
\end{align}
Let $ \mathbb{Z}_{\geq} = \mathbb{Z} \cup \{ 0 \} $. 
Then the state at time $n \in \mathbb{Z}_{\ge}$ and location ${\bf x} \in V( T^d_N )$ can be expressed by a $2d$-dimensional vector:
\begin{align*}
\Psi_{n}( {\bf x} )=
\begin{bmatrix}
\Psi^{1}_{n}( {\bf x} ) \\ \Psi^{2}_{n}( {\bf x} ) \\ \vdots \\ \Psi^{2d}_{n}( {\bf x} ) 
\end{bmatrix} 
\in \mathbb{C}^{2d}.
\end{align*} 
For $\Psi_n : V( T^d_N ) \longrightarrow \mathbb{C}^{2d} \ (n \in \mathbb{Z}_{\geq})$, from (1), the evolution of the walk is defined by 
\begin{align}
\Psi_{n+1}( {\bf x} ) \equiv (M_A \Psi_{n})( {\bf x} )=\sum_{j=1}^{d}\Big(P_{2j-1}A\Psi_{n}( {\bf x} + {\bf e}_j)+P_{2j}A\Psi_{n}( {\bf x} - {\bf e}_j)\Big).
\end{align}

Now, we define the {\em walk-type zeta function} by 
\begin{align}
\overline{\zeta} \left(A, T^d_N, u \right) = \det \Big( I_{2d N^d} - u M_A \Big)^{-1/N^d}.
\label{satosan01}
\end{align}
In general, for a $d_c \times d_c$ coin matrix $A$, we put  
\begin{align*}
\overline{\zeta} \left(A, T^d_N, u \right) = \det \Big(I_{d_c N^d} - u M_A \Big)^{-1/N^d}.
\end{align*}

Komatsu et al. \cite{K2} 
obtained the following result. 
\begin{theorem}[Komatsu et al. \cite{K2}] 
\begin{align*}
\overline{\zeta} \left(A, T^d_N, u \right) ^{-1}
&= \exp \left[ \frac{1}{N^d} \sum_{\widetilde{ {\bf k} } \in \widetilde{\mathbb{K}}_N^d} \log \left\{ \det \Big( F(\widetilde{ {\bf k} }, u) \Big) \right\} \right],
\\
\lim_{N \to \infty} \overline{\zeta} \left(A, T^d_N, u \right) ^{-1}
&=
\exp \left[ \int_{[0,2 \pi)^d} \log \left\{ \det \Big( F \left( \Theta^{(d)}, u \right)  \Big) \right\} d \Theta^{(d)}_{unif} \right],
\end{align*}
where $ \tilde{\mathbb{K}}_N = \{ 0, 2 \pi /N, \ldots , 2 \pi (N-1)/N \} $, $\Theta^{(d)} = (\theta_1, \theta_2, \ldots, \theta_d) (\in [0, 2 \pi)^d)$ 
and $d \Theta^{(d)}_{unif}$ denotes the uniform measure on $[0, 2 \pi)^d$, that is,
\begin{align*}
d \Theta^{(d)}_{unif} = \frac{d \theta_1}{2 \pi } \cdots \frac{d \theta_d}{2 \pi } .  
\end{align*}
Furthermore,  
\begin{align*}
F \left( {\bf w} , u \right) = I_{2d} - u \widehat{M}_A ( {\bf w} ) \ \ \  and \ \ \                           
\widehat{M}_A( {\bf w} )=\sum_{j=1}^{d} \Big( e^{i w_j } P_{2j-1} A + e^{-i w_j } P_{2j} A \Big) , 
\end{align*}
with $ {\bf w} = (w_1, w_2, \ldots, w_d) \in \RM^d$.
\end{theorem}

\section{A lower bound on the exponential rate in the SIS model} 

At first, we state the stochastic SIS model on a digraph. 
Let $D=(V(D),A(D))$ be a connected digraph with $N$ vertices and $M$ arcs, where $V(D)$ and $A(D)$ are the set 
of vertices and arcs of $D$, respectively. 
At any given continuous time $t \geq 0$, each vertex is in one of two possible states, namely, susceptible or infected. 
An infected vertex $v$ stochastically transits to the susceptible state at a constant instantaneous rate of $\delta {}_v >0$, 
which is called the {\em recovery rate} of $v$. 
Whenever $v$ is susceptible and its infected in-neighbor $w$ is infected, then $w$ stochastically and independently 
infects $v$ at a constant instantaneous rate of $\beta {}_{(w,v)} $. 
The value $\beta {}_{(w,v)} $ is called the {\em infection rate}. 

Let $p_v (t)\ (v \in V(D))$ be the probability that $v$ is infected at time $t$. 
Then the {\em decay rate} of the SIS model is defined by 
\[
\gamma =- \limsup_{t \rightarrow \infty} \frac{\log \sum_{v \in V(D)} p_v (t)}{t} ,  
\]
where all vertices are assumed to be infected at $t=0$. 

Next, let $V(D)= \{ v_1 , \ldots , v_N \} $ and $A(D)= \{ e_1 , \ldots , e_M \} $, and 
let an $N \times M$ matrix $C$ be the {\em incidence matrix} of $D$. 
Furthermore, we define the {\em non-backtracking matrix} $H=( H_{ef} )_{e,f \in A(D)} $ of $D$ as follows: 
\[
H_{ef} = 
\left\{
\begin{array}{ll}
1 & \mbox{if $t(e)=o(f)$ and $f \neq e^{-1} $, } \\
0 & \mbox{otherwise. }
\end{array}
\right. 
\]
The matrix $H$ is a digraph version of the edge matrix of a graph.  

A real matrix $A$ is callled {\em nonnegative}, denoted $A \geq {\bf 0} $, if all entries of $A$ are nonnegative. 
For two matrices $A$ and $B$ of same size, we say that $A \leq B$ if $B-A \geq {\bf 0} $. 
A square matrix $A$ is called {\em Metzler} if all its off-diagonal entries are nonnegative \cite{FR}. 
For the Metzler matrix $A$, the maximum real part of the eigenvalues of $A$ is denoted by $\lambda {}_{\max} (A)$ (see \cite{FR}). 
Let \( {\bf M}_{1} \oplus \cdots \oplus {\bf M}_{s} \) be the 
block diagonal sum of square matrices 
${\bf M}_{1}, \cdots , {\bf M}_{s}$. 
If \( {\bf M}_{1} = {\bf M}_{2} = \cdots = {\bf M}_{s} = {\bf M} \),
then we write 
\( s \circ {\bf M} = {\bf M}_{1} \oplus \cdots \oplus {\bf M}_{s} \).
Then an $(N+M) \times (N+M)$ ${\cal A} $ is defined as follows: 
\[ 
{\cal A} = 
\left[ 
\begin{array}{cc}
- {\bf D} & C_+ B^{\prime } \\ 
D^{\prime }_2 C^T_- & H^T B^{\prime } - B^{\prime } - D^{\prime }_1 - D^{\prime }_2 
\end{array} 
\right] 
, 
\] 
where 
\[
C_+ = \max (C, {\bf 0} ) , \ C_- = \max (-C, {\bf 0} ) 
\]
and 
\begin{align*}
{\bf D} & = (\delta {}_{v_1 } ) \otimes \cdots \otimes (\delta {}_{v_N } ) ,  \\ 
D^{\prime }_1 &= (\delta {}_{o(e_1 )} ) \otimes \cdots \otimes (\delta {}_{o(e_M)} ) ,  \\ 
D^{\prime }_2 &= (\delta {}_{t(e_1 )} ) \otimes \cdots \otimes (\delta {}_{t(e_m)} ) ,  \\ 
B^{\prime } &= ( \beta {}_{e_1 } ) \otimes \cdots \otimes ( \beta {}_{e_m} ) . 
\end{align*}
Two matrices $C_+ $ and $C_- $ denote the positive and negative parts of the incidence matrix $C$, respectively. 
Note that ${\cal A} $ is Metzler. 

Masuda, Preciado and Ogura \cite{MPO} presented a new lower bound on the exponential rate by using the leading eigenvalue of ${\cal A} $.

\begin{theorem}[Masuda, Preciado and Ogura \cite{MPO}]  
A lower bound on the decay rate $\gamma $ is given as follows: 
\[
\gamma \geq - \lambda {}_{\max} ( {\cal A} ) . 
\]
\end{theorem}

Furthermore, we use the Weinstein-Aronszajn identity (see \cite{F}).

\begin{theorem}[the Weinstein-Aronszajn identity]
If ${\bf A}$ and ${\bf B}$ are an $r \times s$ and an $s \times r$ 
matrix, respectively, then we have 
\[
\det ( {\bf I}_{r} - {\bf A} {\bf B} )= 
\det ( {\bf I}_s - {\bf B} {\bf A} ) . 
\] 
\end{theorem}

\section{A determinant formula for the matrix ${\cal A} $ of a digraph}

At first, we state the properties for various matrices composing the Metzler matrix ${\cal A} $. 

Let $D$ be a connected graph with $N$ vertices $v_1 , \cdots , v_N $ 
and $M$ arcs, and $A(D)$ the arc set of $D$. 
Furthermore, let two $N \times N$ matrices $B=(B_{uv} )$ and ${\bf D} = {\bf D} (D)=( D_{uv} )$ be 
given as follows: 
\[
B_{uv} =\left\{
\begin{array}{ll}
\beta {}_{uv} & \mbox{if $(u,v) \in A(D)$, } \\
0 & \mbox{otherwise, }
\end{array}
\right.
\ 
D_{uv} =\left\{
\begin{array}{ll}
\delta {}_u & \mbox{if $u=v$, } \\
0 & \mbox{otherwise,}
\end{array}
\right.
\] 
where $\beta{}_{uv} >0, \ \delta {}_u >0$ for each $e=(u,v) \in A(D)$ and $u \in V(D)$. 
Set $\beta{}_e = \beta{}_{uv} $ for $e=(u,v) \in A(D)$. 

The {\em incidence matrix} $C=( C_{ue} )$ of $D$ is the $N \times M$ matrix given as follows: 
\[
C_{eu} =\left\{
\begin{array}{ll}
1 & \mbox{if $t(e)=u$, } \\
-1 & \mbox{if $o(e)=u$. } \\
0 & \mbox{otherwise. }
\end{array}
\right. 
\]
Furthermore, two matrices $C_+ $ and $C_- $ are defined as follows: 
\[
C_+ = \max (C, {\bf 0} ) , \ C_- = \max (-C, {\bf 0} ) . 
\]
Then we have 
\[
(C_+ )_{eu} =\left\{
\begin{array}{ll}
1 & \mbox{if $t(e)=u$, } \\
0 & \mbox{otherwise, }
\end{array}
\right.
\ 
( C_- )_{eu} =\left\{
\begin{array}{ll}
1 & \mbox{if $o(e)=u$, } \\
0 & \mbox{otherwise. }
\end{array}
\right.
\] 
That is, $C_+ $ and $C_- $ are the terminus-arc incidence matrix and the origin-arc incidence matrix 
of $D$, respectively. 
For brevity, we set 
\[
K=C_+ , \ L=C_- . 
\]

Three $M \times M$ diagonal matrices $B^{\prime } =( B^{\prime }_{ef} )$, $D_1 =( D^{(1)}_{ef} )$ and 
$D^{\prime }_2 =( D^{(2)}_{ef} )$ are given by 
\[
B^{\prime }_{ef} =\left\{
\begin{array}{ll}
\beta {}_e & \mbox{if $e=f$, } \\
0 & \mbox{otherwise, }
\end{array}
\right.
\ 
D^{(1)}_{ef} =\left\{
\begin{array}{ll}
\delta {}_{o(e)} & \mbox{if $e=f$, } \\
0 & \mbox{otherwise, }
\end{array}
\right. 
\ 
D^{(2)}_{ef} =\left\{
\begin{array}{ll}
\delta {}_{t(e)} & \mbox{if $e=f$, } \\
0 & \mbox{otherwise. }
\end{array}
\right.
\]
Then $KB^{\prime }$ and $ D^{\prime}_2 L^T$ are the following matrices: 
\[
(KB^{\prime } )_{ue} =\left\{
\begin{array}{ll}
\beta {}_e & \mbox{if $t(e)=u$, } \\
0 & \mbox{otherwise, }
\end{array}
\right.
\ 
( D^{\prime}_2 L^T)_{ue} =\left\{
\begin{array}{ll}
\delta {}_{t(e)} & \mbox{if $o(e)=u$, } \\
0 & \mbox{otherwise. }
\end{array}
\right. 
\]

Next, let 
\[
M_0 =| \{ e \in A(D) \mid  e^{-1} \notin A(D) \} | , \ 
M_1 =| \{ e \in A(D) \mid  e^{-1} \in A(D) \} |/2 . 
\]
Note that $M=M_0 +2 M_1 $. 
Furthermore, let 
\[
A(D)= \{ e_1 ,\ldots , e_{M_0} , f_1 , \ldots , f_{M_1} ,  f^{-1}_1 , \ldots , f^{-1}_{M_1} \} . 
\] 
We consider the matrices under this order. 
The matrix $M \times M$ matrix $J$ is defined as follows: 
\[
J=  
\left[ 
\begin{array}{ccc}
{\bf 0} & {\bf 0} & {\bf 0} \\ 
{\bf 0} & {\bf 0} & {\bf I}_M \\  
{\bf 0} & {\bf I}_M & {\bf 0}  
\end{array} 
\right] 
. 
\]

Then the following result holds.

\newtheorem{proposition}{Proposition} 
\begin{proposition} 
The non-backtracking matrix $H$ of $D$ is expressed as follows: 
\[
H= K^T L -J . 
\]
\end{proposition}

Here, we set $F= K^T L$. 
Then we have 
\[
H=F-J, \ H^T = L K^T -J =F^T -J . 
\]
The non-backtracking matrix of a digraph is a digraph version of the edge matrix of a graph. 
Furthermore, the Metzler matrix ${\cal A} = {\cal A} (D)$ is written as follows: 
\[
{\cal A} =  
\left[ 
\begin{array}{cc}
- {\bf D} &  KB^{\prime } \\ 
D^{\prime}_2 L^T & H^T B^{\prime } - B^{\prime } - D^{\prime }_1 - D^{\prime }_2   
\end{array} 
\right] 
. 
\]
Set 
\[
E= B^{\prime } + D^{\prime }_1 + D^{\prime }_2. 
\]

Then the following result holds for $\det ( I_{N+M} -u {\cal A} )$.
The {\em Kronecker product} $ {\bf A} \bigotimes {\bf B} $
of matrices {\bf A} and {\bf B} is considered as the matrix 
{\bf A} having the element $a_{ij}$ replaced by the matrix $a_{ij} {\bf B}$.

\begin{theorem}
Let $D$ be a connected digraph with $N$ vertices and $M$ arcs. 
Furthermore, let 
\[
A(D)= \{ e_1 ,\ldots , e_{M_0} , f_1 , \ldots , f_{M_1} ,  f^{-1}_1 , \ldots , f^{-1}_{M_1} \} ,  
\]  
where $e^{-1}_i \notin A(D) \ ( \leq i \leq M_0 )$. 
Then  
\[ 
\det (I_{n+M} -u {\cal A} )= \prod^{M_0}_{j=1} (1+ \gamma {}_{e_j} u) \prod^{M_1}_{k=1} G_k 
\det ( I_N +uD-u A^T_{nsym} -u A_{sym} + u^2 \tilde{D} ) ,  
\] 
where $\gamma {}_e = \beta {}_e + \delta {}_{o(e)} + \delta {}_{t(e)} $ for each $e \in A(D)$, and 
\[
G_k =1+( \gamma {}_{f_k} + \gamma {}_{f^{-1}_k} )u + u^2 ( \gamma {}_{f_k} \gamma {}_{f^{-1}_k} - \beta {}_{f_k} \beta {}_{f^{-1}_k } ) 
\ (1 \leq k \leq M_1 ) . 
\]
\end{theorem}

{\em Proof }.  At first, we have 
\begin{align*}
& \det (  I_{N+M} -u {\cal A} ) \\ 
&= \det 
\left[ 
\begin{array}{cc}
I_N +u {\bf D} & -u KB^{\prime } \\ 
-u D^{\prime}_2 L^T & I_M -u( H^T B^{\prime } -E)  
\end{array} 
\right] 
\\
& \ \ \ \ \ \times \det 
\left[ 
\begin{array}{cc}
I_N & u(I_N +u {\bf D} )^{-1} KB^{\prime } \\ 
{\bf 0} & I_M  
\end{array} 
\right] 
\\
&= \det 
\left[ 
\begin{array}{cc}
I_N +u {\bf D} & {\bf 0} \\ 
-u D^{\prime}_2 L^T & I_M -u( H^T B^{\prime } -E)- u^2  D^{\prime}_2 L^T (I_N +u {\bf D} )^{-1} KB^{\prime }   
\end{array} 
\right] 
. 
\end{align*}
By Proposition 1, we get  
\begin{align*}
& \det (  I_{N+M} -u {\cal A} ) \\ 
&= \det (I_N +u {\bf D} ) \det (I_M -u( H^T B^{\prime } -E)- u^2  D^{\prime}_2 L^T (I_N +u {\bf D} )^{-1} KB^{\prime } ) \\ 
&= \det (I_N +u {\bf D} ) \det (I_M +uE-u(L^T K B^{\prime } -J B^{\prime } )- u^2  D^{\prime}_2 L^T (I_N +u {\bf D} )^{-1} KB^{\prime } ) \\ 
&= \det (I_N +u {\bf D} ) \det (I_M +u(J B^{\prime } +E)-u(L^T +u  D^{\prime}_2 L^T (I_N +u {\bf D} )^{-1} )KB^{\prime } ) \\ 
&= \det (I_N +u {\bf D} ) \det (I_M -u(L^T +u D^{\prime}_2 L^T (I_N +u {\bf D} )^{-1} )KB^{\prime } (I_M +u(J B^{\prime } +E)^{-1} ) \\ 
& \ \ \ \ \ \times \det ( I_M +u(J B^{\prime }+E)) .  
\end{align*}

From Theorem 4, we obtain  
\begin{align*}
& \det (  I_{N+M} -u {\cal A} ) \\ 
&= \det (I_N +u {\bf D} ) \det (I_N -uKB^{\prime } (I_M +u(J B^{\prime }+E))^{-1} (L^T +u D^{\prime}_2 L^T (I_N +u {\bf D} )^{-1} )) \\ 
& \ \ \ \ \ \times \det ( I_M +u(J B^{\prime } +E)) \\ 
&= \det ((I_N +u {\bf D} )-uKB^{\prime } (I_M +u(J B^{\prime }+E))^{-1} (L^T (I_N +u {\bf D} )+u D^{\prime}_2 L^T )) \\ 
& \ \ \ \ \ \times \det ( I_M +u(J B^{\prime } +E)) \\ 
&= \det (I_N +u {\bf D} -uKB^{\prime } (I_M +u(J B^{\prime }+E))^{-1} (L^T +u(L^T {\bf D} + D^{\prime}_2 L^T )) \\ 
& \ \ \ \ \ \times \det ( I_M +u(J B^{\prime } +E)) . 
\end{align*}

Next, we consider $ \det ( I_M +u(J B^{\prime } +E))$. 
Let three diagonal matrices $B_0 $, $B_1 $ and $B_2 $ be given as follows: 
\[ 
( B_0 )_{e_i e_i} = \beta {}_{e_i } , \ ( B_1 )_{f_j f_j} = \beta {}_{f_j }, \ ( B_2 )_{f^{-1}_j f^{-1}_j } = \beta {}_{f^{-1}_j } 
\ (1 \leq i \leq M_0 ; 1 \leq j \leq M_1 ) . 
\] 
Then we get  
\[
(I_{N+M} +J)B^{\prime } =  
\left[ 
\begin{array}{ccc}
B_0 & {\bf 0} & {\bf 0} \\ 
{\bf 0} & B_1 & B_2 \\  
{\bf 0} & B_1 & B_2  
\end{array} 
\right] 
. 
\]
Thus, we have 
\[
JB^{\prime } +E=  
\left[ 
\begin{array}{ccc}
F_0 & {\bf 0} & {\bf 0} \\ 
{\bf 0} & F_1 & B_2 \\  
{\bf 0} & B_1 & F_2  
\end{array} 
\right] 
,  
\]
where $F_0 $, $F_1 $ and $F_2 $ are $M_ 0 \times M_0 $, $M_1 \times M_1 $ and $M_1 \times M_1 $ diagonal matrices 
such that 
\[
( F_0 )_{e_i e_i} = \gamma {}_{e_i } , \ ( F_1 )_{f_j f_j} = \gamma {}_{f_j }, \ ( F_2 )_{f^{-1}_j f^{-1}_j } = \gamma {}_{f^{-1}_j } 
\ (1 \leq i \leq M_0 , 1 \leq j \leq M_1 ) . 
\] 
Therefore, 
\begin{align*}
& \det ( I+u( JB^{\prime } +E)) \\ 
&= \det 
\left[ 
\begin{array}{ccc}
I+u F_0 & {\bf 0} & {\bf 0} \\ 
{\bf 0} & I+u F_1 & u B_2 \\  
{\bf 0} & u B_1 & I+u F_2  
\end{array} 
\right] 
\\
&= \det (I+u F_0 )  \det 
\left[ 
\begin{array}{cc}
I+u F_1 & u B_2 \\  
u B_1 & I+u F_2  
\end{array} 
\right] 
\det 
\left[ 
\begin{array}{cc}
I & -u(I+u F_1 )^{-1} B_2 \\  
{\bf 0} & I  
\end{array} 
\right] 
\\
&= \det (I+u F_0 ) \det 
\left[ 
\begin{array}{cc}
I+u F_1 & {\bf 0} \\  
u B_1 & I+u F_2 -u^2 B_1 (I+u F_1 )^{-1} B_2   
\end{array} 
\right] 
\\
&= \det (I+u F_0 ) \det (I+u F_1 ) \det (I+u F_2 -u^2 B_1 (I+u F_1 )^{-1} B_2 ) . 
\end{align*}

But, we get  
\[
B_1 (i+u F_1 )^{-1} B_2 = \oplus^{M_1}_{j=1}  \left( \frac{ \beta {}_{f_j} \beta {}_{f^{-1}_j } }{1+ \gamma {}_{f_j } u} \right) . 
\]
Thus, we have 
\begin{align*}
&  \det ( I+u( JB^{\prime } +E)) \\ 
&= \det (I+u F_0 ) \det ((I+u F_1 )(I+u F_2 ) -u^2 (I+u F_1 ) B_1 (I+u F_1 )^{-1} B_2 ) . 
\end{align*}
Since   
\[
(I+u F_1 ) B_1 (I+u F_1 )^{-1} B_2 = \oplus^{M_1}_{j=1}  \left( \beta {}_{f_j} \beta {}_{f^{-1}_j } \right) ,  
\]
\begin{align*}
&  \det ( I+u( JB^{\prime } +E)) \\ 
&= \det (I+u F_0 ) \det ((I+u F_1 )(I+u F_2 ) -u^2 \oplus^{M_1}_{j=1}  ( \beta {}_{f_j} \beta {}_{f^{-1}_j } )) \\ 
&= \prod^{M_0}_{j=1} (1+ \gamma {}_{e_j} u) \prod^{M_1}_{k=1} (1+ \gamma {}_{f_k} u)(1+ \gamma {}_{f^{-1}_k} u) 
- \beta {}_{f_k} \beta {}_{f^{-1}_k } ) \\
&= \prod^{M_0}_{j=1} (1+ \gamma {}_{e_j} u) \prod^{M_1}_{k=1} (1+( \gamma {}_{f_k} + \gamma {}_{f^{-1}_k} )u 
+ u^2 ( \gamma {}_{f_k} \gamma {}_{f^{-1}_k} - \beta {}_{f_k} \beta {}_{f^{-1}_k } )) . 
\end{align*}

Now, we consider $( I_M +u(J B^{\prime } +E))^{-1} $ by the Schur complement on the inverse matrix (see \cite{Z}). 
If the suitable matrices are nonsingular, then the inverse matrix $U^{-1} $ of a square matrix 
\[
U= 
\left[ 
\begin{array}{cc}
A & B \\  
C & D 
\end{array} 
\right] 
, 
\]
\[
U^{-1} = 
\left[ 
\begin{array}{cc}
(A-B D^{-1} C)^{-1} &  -A^{-1} B(D-CA^{-1} B)^{-1} \\ 
 -(D-CA^{-1} B)^{-1} CA^{-1} & (D-C A^{-1} B)^{-1}
\end{array} 
\right] 
.  
\]
Substituting $A=I+ F_1 u$, $B=u B_2 $, $C=u B_1 $ and $D=I+F_2 u$, we have 
\[
(A-B D^{-1} C)^{-1} = \oplus^{M_1}_{j=1} \left( \frac{1+ \gamma {}_{ f^{-1}_k } u}{G_k } \right) , 
\]
\[
(D-C A^{-1} B)^{-1} = \oplus^{M_1}_{j=1} \left( \frac{1+ \gamma {}_{ f_k } u}{G_k } \right) , 
\]  
\[
A^{-1} B(D-CA^{-1} B)^{-1} = \oplus^{M_1}_{j=1} \left( \frac{u \beta {}_{ f^{-1}_k } }{G_k } \right) , 
\]
\[
(D-CA^{-1} B)^{-1} CA^{-1} = \oplus^{M_1}_{j=1} \left( \frac{u \beta {}_{ f_k } }{G_k } \right) , 
\]
that is, 
\[
\left[ 
\begin{array}{cc}
I+F_1 u & uB_2 \\  
uB_1 & I+F_2 u 
\end{array} 
\right]^{-1}  
= 
\left[ 
\begin{array}{cc}
\oplus^{M_1}_{j=1} \left( \frac{1+ \gamma {}_{ f^{-1}_k } u}{G_k } \right) & \oplus^{M_1}_{j=1} \left( \frac{u \beta {}_{ f^{-1}_k } }{G_k } \right) \\  
\oplus^{M_1}_{j=1} \left( \frac{u \beta {}_{ f_k } }{G_k } \right) & \oplus^{M_1}_{j=1} \left( \frac{1+ \gamma {}_{ f_k } u}{G_k } \right)  
\end{array} 
\right]
. 
\]
Thus, 
\[ 
\begin{array}{rcl}
\  &  & ( I_M +u(J B^{\prime }+E))^{-1} \\ 
\  &   &                \\ \  &   &                \\ 
\  & = & 
\left[ 
\begin{array}{ccc}
\oplus^{M_0}_{j=1} \left( \frac{1}{1+ \gamma {}_{ e^{-1}_j } u} \right) & {\bf 0} & {\bf 0} \\  
{\bf 0} & \oplus^{M_1}_{j=1} \left( \frac{1+ \gamma {}_{ f^{-1}_k } u}{G_k } \right) & \oplus^{M_1}_{j=1} \left( \frac{u \beta {}_{ f^{-1}_k } }{G_k } \right) \\  
{\bf 0} & \oplus^{M_1}_{j=1} \left( \frac{u \beta {}_{ f_k } }{G_k } \right) & \oplus^{M_1}_{j=1} \left( \frac{1+ \gamma {}_{ f_k } u}{G_k } \right)  
\end{array} 
\right] 
. 
\end{array}
\]

Now, let 
\[
M= KB^{\prime } (I_M +u(J B^{\prime } +E))^{-1} (L^T +u(L^T {\bf D} + D^{\prime}_2 L^T )) .
\]
Furthermore, put  
\[
\epsilon {}_e = \delta {}_{o(e)} + \delta {}_{t(e)} \ for \ e \in A(D) . 
\]
Then we get  
\[
(L^T +u(L^T {\bf D} + D^{\prime}_2 L^T ))_{eu} = 
\left\{
\begin{array}{ll}
1+u \epsilon {}_e & \mbox{if $o(e)=u$, } \\
0 & \mbox{otherwise. }
\end{array}
\right. 
\]
Thus, the entries of $M$ are given as follows: 
\[
M_{uv} = 
\left\{
\begin{array}{ll}
\beta {}_{e^{-1}} (1+ \epsilon {}_{e^{-1}} u)/(1+ \gamma {}_{e^{-1}} u) & \mbox{if $e=(u,v) \in A(D)$ and $e^{-1} \notin A(D)$, } \\
\beta {}_{e^{-1}} (1+ \gamma {}_e u)(1+ \epsilon {}_{e^{-1}} u)/ G_e & \mbox{if $e=(u,v) , e^{-1}  \in A(D)$, } \\
-u \sum_{o(e)=u} \beta {}_e \beta {}_{e~{-1}} (1+ \epsilon {}_e u)/ G_e & \mbox{if $u=v$, } \\
0 & \mbox{otherwise. }
\end{array}
\right. 
\] 
Therefore, it follows that 
\[
\det ( I_{N+M} -u {\cal A} )= \prod^{M_0}_{j=1} (1+ \gamma {}_{e_j} u) \prod^{M_1}_{k=1} G_k 
\det ( I_N +uD-uM) . 
\]

Now, we define three $N \times N$ matrices $A_{nsym} $, $A_{sym} $ and $\tilde{D} $ as follows: 
\begin{align*}
(A_{nsym} )_{uv} &= 
\left\{
\begin{array}{ll}
( \beta {}_{e~{-1}} (1+ \epsilon {}_{e^{-1}} u)/(1+ \gamma {}_{e^{-1}} u) & \mbox{if $e=(u,v)  \in A(D)$ and $e^{-1} \notin A(D)$, } \\
0 & \mbox{otherwise, }
\end{array}
\right. 
\\ 
(A_{sym} )_{uv} &=
\left\{
\begin{array}{ll}
( \beta {}_{e~{-1}} (1+ \gamma {}_e u)(1+ \epsilon {}_{e^{-1}} u)/ G_e & \mbox{if $e=(u,v) , e^{-1}  \in A(D)$, } \\
0 & \mbox{otherwise, }
\end{array}
\right. 
\\ 
(\tilde{D} )_{uv} &= 
\left\{
\begin{array}{ll}
-u \sum_{o(e)=u} ( \beta {}_e \beta {}_{e~{-1}} (1+ \epsilon {}_e u)/ G_e & \mbox{if $u=v$, } \\
0 & \mbox{otherwise. }
\end{array}
\right. 
\end{align*}
Hence,  
\[ 
\det (I_{n+M} -u {\cal A} )= \prod^{M_0}_{j=1} (1+ \gamma {}_{e_j} u) \prod^{M_1}_{k=1} G_k 
\det ( I_N +uD-u A_{nsym} -u A^T_{sym} + u^2 \tilde{D} ) . 
\]
$\Box$

\section{The spectrum for the Metzler matrix of the symmetric digraph of a graph} 

Let $D$ be a connected digraph with $N$ vertices and $M$ arcs, and let $\beta {}_e = \beta $, \ $\delta {}_u = \delta $ 
for each $e \in A(D)$ and $u \in V(D)$. 
Then we have 
\[
\gamma {}_e = \beta +2 \delta - \gamma \ (const), \ \epsilon {}_e = 2 \delta = \epsilon \ (const) \ 
for \ each \ e \in A(D)  
\]
and 
\[
G_j =(1+( \beta +2 \delta )u)^2 - \beta {}^2 u^2 =1+2( \beta +2 \delta )u+4 \delta ( \delta +1) u^2 =G 
\ \ \ (1 \leq j \leq M_1 ) . 
\]
Thus, we have 
\[
A_{nsym} = \frac{ \beta (1+ \epsilon u)}{1+ \gamma u} A_0 , \ \ \ 
A_{sym} = \frac{ \beta (1+ \epsilon u)(1+ \gamma u)}{G} A_1 ,  
\]
\[
\tilde{D} = \frac{ \beta {}^2 (1+ \epsilon u)}{G} D_0 , 
\]
where three matrices $A_0 $, $A_1 $ and $D_0 $ are given as follows: 
\[
(A_{0} )_{uv} = 
\left\{
\begin{array}{ll}
1 & \mbox{if $(u,v)  \in A(D)$ and $(v,u) \notin A(D)$, } \\
0 & \mbox{otherwise, }
\end{array}
\right. 
\ 
(A_1 )_{uv} =
\left\{
\begin{array}{ll}
1 & \mbox{if $(u,v) , (v,u) \in A(D)$, } \\
0 & \mbox{otherwise, }
\end{array}
\right. 
\]
\[
( D_0 )_{uv} = 
\left\{
\begin{array}{ll}
\deg {}^+ \ u & \mbox{if $u=v$, } \\
0 & \mbox{otherwise. }
\end{array}
\right. 
\]
Here, $\deg {}^+ \ u =| \{ e \in A(D) \mid o(e)=u \} |$. 
Thus, we have 
\begin{align*}
& \det (I_{n+M} -u {\cal A} )=(1+ \gamma u)^{M_0} G^{M_1} 
\\ 
& \ \ \ \times \det ( I_N + \delta u I_N -u \frac{ \beta (1+ \epsilon u)}{1+ \gamma u} A^T_0 
-u \frac{ \beta (1+ \gamma u)(1+ \epsilon u)}{G} A_1 +u^2 \frac{ \beta {}^2 (1+ \epsilon u)}{G} D_0 ) 
\\ 
&=(1+ \gamma u)^{M_0 -N} G^{M_1 -N} 
\\
& \ \ \ \times \det ((1+ \delta u)(1+ \gamma u)(1+2( \beta +2 \delta )u+4 \delta ( \delta +1) u^2 ) I_N 
\\
& \ \ \ - \beta u (1+ \epsilon u)(1+2( \beta +2 \delta )u+4 \delta ( \delta +1) u^2 ) A^T_0 -u \beta (1+ \gamma u)^2 (1+ \epsilon u) A_1 
+ \beta {}^2 (1+ \gamma u)(1+ \epsilon u) D_0 ) 
\\
&=(1+ \gamma u)^{M_0 -N} (1+2( \beta +2 \delta )u+4 \delta ( \delta +1) u^2 )^{M_1 -N} 
\\
& \ \ \ \times \det ((1+ \delta u)(1+( \beta +2 \delta )u)(1+2( \beta +2 \delta )u+4 \delta ( \delta +1) u^2 ) I_N 
\\
& \ \ \ - (1+2 \delta u) \beta u ((1+2( \beta +2 \delta )u+4 \delta ( \delta +1) u^2 ) A^T_0 +(1+( \beta +2 \delta )u)^2 A_1 )
- \beta (1+( \beta +2 \delta ) u) D_0 ) . 
\end{align*}

Therefore, we obtain

\newtheorem{corollary}{Corollary}
\begin{corollary}
Let $D$ be a connected graph with $N$ vertices $v_1 , \cdots , v_N $ 
and $M$ arcs, and let $\beta {}_e = \beta $, \ $\delta {}_u = \delta $ for each $e \in A(D)$ and $u \in V(D)$. 
Then 
\begin{align*}
& \det (I_{n+M} -u {\cal A} )=(1+ \gamma u)^{M_0 -N} (1+2( \beta +2 \delta )u+4 \delta ( \delta +1) u^2 )^{M_1 -N} 
\\
& \ \ \ \times \det ((1+ \delta u)(1+( \beta +2 \delta )u)(1+2( \beta +2 \delta )u+4 \delta ( \delta +1) u^2 ) I_N 
\\
& \ \ \ - (1+2 \delta u) \beta u ((1+2( \beta +2 \delta )u+4 \delta ( \delta +1) u^2 ) A^T_0 +(1+( \beta +2 \delta )u)^2 A_1 ) 
- \beta  (1+( \beta +2 \delta ) u) D_0 ) . 
\end{align*}
\end{corollary}

Next, let $D=G$ be a connected (undirected) graph with $N$ vertices and $M$ edges, 
and let $\beta {}_e = \beta $, \ $\delta {}_u = \delta $ for each $e \in D(G)$ and $u \in V(G)$. 
Then we have 
\[
M_0 =0, \ M_1 =M, \ |D(G)|=2M . 
\]
Furthermore, we get  
\[
A_0 = {\bf 0}, \ A_1 ={\bf A} (G) \ and \ D_0 = {\bf D}_G . 
\]

Thus,

\begin{corollary}
Let $D=G$ be a connected (undirected) graph with $N$ vertices and $M$ edges, 
and let $\beta {}_e = \beta $, \ $\delta {}_u = \delta $ for each $e \in D(G)$ and $u \in V(G)$. 
Then 
\begin{align*}
& \det (I_{n+2M} -u {\cal A} )=(1+2( \beta +2 \delta )u+4 \delta ( \delta +1) u^2 )^{M-N} 
\\
&\times \det ((1+ \delta u)(1+2( \beta +2 \delta )u+4 \delta ( \delta +1) u^2 ) I_N 
\\
&- (1+2 \delta u) \beta u (1+( \beta +2 \delta )u) {\bf A} (G) + \beta {}^2 (1+2 \delta u) u^2 {\bf D}_G ) . 
\end{align*}
\end{corollary}

Finally, let $G$ be a connected $d$-regular graph with $N$ vertices and $M$ edges, 
and let $\beta {}_e = \beta $, \ $\delta {}_u = \delta $ for each $e \in D(G)$ and $u \in V(G)$. 
By Corollary 2, we have 
\begin{align*}
& \det (I_{n+2M} -u {\cal A} )=(1+2( \beta +2 \delta )u+4 \delta ( \delta +1) u^2 )^{M-N}  
\\ 
& \ \ \ \times \det ( (1+ \delta u)(1+2( \beta +2 \delta )u+4 \delta ( \delta +1) u^2 ) I_N 
\\
& \ \ \ - (1+2 \delta u) \beta u (1+( \beta +2 \delta )u) {\bf A} (G) +d \beta {}^2 (1+2 \delta u) u^2 I_N ) 
\\ 
&= (1+2( \beta +2 \delta )u+4 \delta ( \delta +1) u^2 )^{M-N}  
\\ 
& \ \ \ \times \det ( \{ 1+(2 \beta +5 \delta )u+(2 \delta ( \beta +4 \delta +2)+d \beta {}^2 ) u^2 
+2 \delta (2 \delta ( \delta +1)+d \beta {}^2 ) u^3 \} I_N 
\\
& \ \ \ - \beta (1+2 \delta u)(1+( \beta +2 \delta )u) {\bf A} (G) ) . 
\end{align*}

Therefore,

\begin{corollary}
Let $D=G$ be a a connected $d$-regular graph with $N$ vertices and $M$ edges, 
and let $\beta {}_e = \beta $, \ $\delta {}_u = \delta $ for each $e \in D(G)$ and $u \in V(G)$. 
Then 
\begin{align*}
& \det (I_{n+2M} -u {\cal A} )=(1+2( \beta +2 \delta )u+4 \delta ( \delta +1) u^2 )^{M-N}  
\\ 
& \ \ \ \times \det ( \{ 1+(2 \beta +5 \delta )u+(2 \delta ( \beta +4 \delta +2)+d \beta {}^2 ) u^2 
+2 \delta (2 \delta ( \delta +1)+d \beta {}^2 ) u^3 \} I_N 
\\
& \ \ \ - \beta (1+2 \delta u)(1+( \beta +2 \delta )u) {\bf A} (G) ) . 
\end{align*}
\end{corollary}

From Corollary 3, the following result holds.

\begin{corollary}
Let $D=G$ be a a connected $d$-regular graph with $N$ vertices and $M$ edges, 
and let $\beta {}_e = \beta $, \ $\delta {}_u = \delta $ for each $e \in D(G)$ and $u \in V(G)$. 
Then 
\begin{align*}
& \det ( \lambda I_{n+2M} - {\cal A} )=( \lambda {}^2 +2( \beta +2 \delta ) \lambda +4 \delta ( \delta +1) )^{M-N}  
\\ 
& \ \ \ \times \det ( \{ \lambda {}^3 +(2 \beta +5 \delta ) \lambda {}^2 +(2 \delta ( \beta +4 \delta +2)+d \beta {}^2 ) \lambda  
+2 \delta (2 \delta ( \delta +1)+d \beta {}^2 ) \} I_N 
\\
& \ \ \ - \beta ( \lambda +2 \delta )( \lambda +( \beta +2 \delta )) {\bf A} (G) ) . 
\end{align*}
\end{corollary}

Substituting $u=1 /\lambda $ to Corollary 4, the spectra for ${\cal A} $ of a regular graph are obtained.

\begin{corollary}
Let $D=G$ be a a connected $d$-regular graph with $N$ vertices and $M$ edges, 
and let $\beta {}_e = \beta $, \ $\delta {}_u = \delta $ for each $e \in D(G)$ and $u \in V(G)$. 
Then the eigenvalues of ${\cal A} $ are given as follows: 
\begin{enumerate} 
\item $3N$ eigenvalues: solutions of 
\begin{align*}
& \lambda {}^3 +(2 \beta - \beta \mu +5 \delta ) \lambda {}^2 +(2 \delta ( \beta +4 \delta +2) +d \beta {}^2 
- \beta ( \beta +4 \delta ) \mu ) \lambda \\
&+ 2 \delta ( 2 \delta ( \delta +1)+d \beta ^2 - \beta ( \beta +2 \delta ) \mu )=0 , \ \ \mu \in {\rm Spec} ({\bf A} (G)) , 
\end{align*}
\item $2(M-N)$ eigenvalues: $\lambda =-( \beta +2 \delta ) \pm \sqrt{ \beta {}^2 +4 \delta ( \beta -1)} $.  
\end{enumerate} 
\end{corollary}

{\bf Proof }.  By Corollary 4, we have 
\begin{align*}
& \det ( \lambda I_{n+2M} - {\cal A} )=( \lambda {}^2 +2( \beta +2 \delta ) \lambda +4 \delta ( \delta +1) )^{M-N}  
\\ 
& \ \ \ \times \prod_{ \mu \in Spec( {\bf A} (G))} \{ \lambda {}^3 +(2 \beta - \beta \mu +5 \delta ) \lambda {}^2 
+(2 \delta ( \beta +4 \delta )+2) +d \beta {}^2  
\\
& \ \ \ - \beta ( \beta +4 \delta ) \mu ) \lambda +2 \delta ( 2 \delta ( \delta +1)+d \beta ^2 - \beta ( \beta +2 \delta ) \mu ) \} . 
\end{align*}

Solving $  \lambda {}^2 +2( \beta +2 \delta ) \lambda +4 \delta ( \delta +1)=0$, we obtain 
\[
\lambda =-( \beta +2 \delta ) \pm \sqrt{ \beta {}^2 +4 \delta ( \beta -1)} . 
\] 
The result follows. 
$\Box$

\section{A walk-type zeta function with respect to the Metzler matrix of the symmetric digraph of a finite torus}

As an application, we present the Walk/Zeta Correspondence (see \cite{K2}) 
for the Metzler matrix of the symmetric digraph of the $d$-dimensional torus $G= T^d_N $. 
For a natural number $N \geq 2$, let $G= T^d_N $ be the $d$-dimensional finite torus (graph). 
Furthermore, assume that $\beta {}_e = \beta $, \ $\delta {}_u = \delta $ for each $e \in D( T^d_N )$ and $u \in V( T^d_N )$. 
Then we have 
\[
n=|V( T^d_N )|= N^d , \ m=|E( T^d_N )|=d N^d . 
\]
Note that $T^d_N $ is a $2d$-regular graph.

Similarly to the Walk/Zeta Correspondence of \cite{K2}, 
we introduce a zeta function for the Metzler matrix ${\cal A} $ of the symmetric digraph of $G= T^d_N $ as follows: 
\[
\overline{\zeta } ( {\cal A} , T^d_N , u)= \det ( {\bf I}_{(2d+1) N^d} -u {\cal A} )^{-1/ N^d } . 
\]
Let ${\bf A} = {\bf A} (G)= {\bf A} ( T^d_N )$. 
It is known that 
\[
{\rm Spec} ( T^d_N )= \Big\{ 2 \sum^d_{j=1} \cos \left( \frac{ 2 \pi k_j }{N} \right) \mid k_j =0,1, \ldots , N-1 \ (1 \leq j \leq d) \Big\} . 
\]
By Theorem 5 and Corollary 3, we have 
\begin{align*}
& \overline{\zeta } ( {\cal A} , T^d_N , u)^{-1} = \det ( {\bf I}_{(2d+1) N^d} -u {\cal A} )^{1/ N^d } \\ 
&= \{ (1+2( \beta +2 \delta )u+4 \delta ( \delta +1) u^2 )^{(2d-1) N^d } 
\times \det ( \{ 1+(2 \beta +5 \delta )u+(2 \delta ( \beta +4 \delta +2)+2d \beta {}^2 ) u^2 \\
& \ \ \ + 2 \delta (2 \delta ( \delta +1)+2d \beta {}^2 ) u^3 \} I_{N^d } - \beta (1+2 \delta u)(1+( \beta +2 \delta )u) {\bf A} (T^d_N ) ) \} {}^{1/N^d } \\ 
&= (1+2( \beta +2 \delta )u+4 \delta ( \delta +1) u^2 )^{2d-1} \prod^{N-1}_{k_1 =0} \cdots \prod^{N-1}_{k_d =0} \{ \{ 1+(2 \beta +5 \delta )u \\
& \ \ \ + 2( \delta ( \beta +4 \delta +2)+d \beta {}^2 ) u^2 +4 \delta ( \delta +1)+d \beta {}^2 ) u^3 \} \\ 
& \ \ \ - 2 \beta (1+2 \delta u)(1+( \beta +2 \delta )u) \sum^d_{j=1} \cos \left( \frac{ 2 \pi k_j }{N} \right) \} {}^{1/N^d } \\ 
&= (1+2( \beta +2 \delta )u+4 \delta ( \delta +1) u^2 )^{2d-1} \\
& \ \ \ \times \exp ( \frac{1}{N^d} \sum^{N-1}_{k_1 =0} \cdots \sum^{N-1}_{k_d =0} \log \{ \{ 1+(2 \beta +5 \delta )u+2( \delta ( \beta +4 \delta +2)+d \beta {}^2 ) u^2 
+4 \delta ( \delta +1)+d \beta {}^2 ) u^3 \} \\ 
& \ \ \ - 2 \beta (1+2 \delta u)(1+( \beta +2 \delta )u)  \sum^d_{j=1} \cos \left( \frac{ 2 \pi k_j }{N} \right) \} ) . 
\end{align*}

When $N \rightarrow \infty$, we get  
\begin{align*}
& \lim_{N \rightarrow \infty} \overline{\zeta } ( {\cal A} , T^d_N , u)^{-1}\\ 
&= (1+2( \beta +2 \delta )u+4 \delta ( \delta +1) u^2 )^{2d-1} \\
& \ \ \ \times \exp  ( \int^{2 \pi }_{0} \cdots \int^{2 \pi }_{0} \log \{ \{ 1+(2 \beta +5 \delta )u+2( \delta ( \beta +4 \delta +2)+d \beta {}^2 ) u^2 
+4 \delta ( \delta +1)+d \beta {}^2 ) u^3 \} \\ 
& \ \ \ - 2 \beta (1+2 \delta u)(1+( \beta +2 \delta )u) \sum^d_{j=1} \cos \theta {}_j \frac{d \theta {}_1 }{ 2 \pi } 
\cdots \frac{d \theta {}_d }{ 2 \pi } \} ) . 
\end{align*}

Therefore, we obtain

\begin{theorem} 
Let $G= T^d_N $ be the d-dimensional finite torus. 
Then 
\begin{align*}
& \overline{\zeta } ( {\cal A} , T^d_N , u)^{-1} \\
&= (1+2( \beta +2 \delta )u+4 \delta ( \delta +1) u^2 )^{2d-1} \\
&\ \ \ \times \exp  \Big[ \frac{1}{N^d} \sum^{N-1}_{k_1 =0} \cdots \sum^{N-1}_{k_d =0} \log \Big\{ \{ 1+(2 \beta +5 \delta )u+2( \delta 
( \beta +4 \delta +2)+d \beta {}^2 ) u^2 \\ 
&\ \ \ + 4 \delta ( \delta +1)+d \beta {}^2 ) u^3 \} I_{N^d } - 2 \beta (1+2 \delta u)(1+( \beta +2 \delta )u)  
\sum^d_{j=1} \cos \left( \frac{ 2 \pi k_j }{N} \right) \Big\}  \Big]   
\end{align*}
and  
\begin{align*}
& \lim_{N \rightarrow \infty} \overline{\zeta } ( {\cal A} , T^d_N , u)^{-1}\\ 
&= (1+2( \beta +2 \delta )u+4 \delta ( \delta +1) u^2 )^{2d-1} \\
& \ \ \ \times \exp  \Big[ \int^{2 \pi }_{0} \cdots \int^{2 \pi }_{0} \log \Big\{ \{ 1+(2 \beta +5 \delta )u+2( \delta ( \beta +4 \delta +2)+d \beta {}^2 ) u^2 
+4 \delta ( \delta +1)+d \beta {}^2 ) u^3 \} \\ 
& \ \ \ - 2 \beta (1+2 \delta u)(1+( \beta +2 \delta )u) \sum^d_{j=1} \cos \theta {}_j \frac{d \theta {}_1 }{ 2 \pi } 
\cdots \frac{d \theta {}_d }{ 2 \pi } \Big\} \Big] . 
\end{align*}
\end{theorem}

\end{document}